\documentclass[a4paper]{article}

\usepackage{amsfonts}
\usepackage{amsmath}
\usepackage{amssymb}
\usepackage{mathrsfs}
\usepackage{amscd}
\usepackage{mathtools}
\usepackage{here}
\usepackage{tikz}
\usepackage[all]{xy}
\usepackage{xcolor}
\usepackage{bbm}
\usepackage{colortbl}
\usepackage{selinput}
\usepackage{hyperref}
\hypersetup{
    colorlinks=true,
    linkcolor=blue,
    filecolor=magenta,
    urlcolor=cyan,
    pdftitle={Overleaf Example},
    pdfpagemode=FullScreen,
    }
\newcommand{\z}{\mathbb{Z}}

\newcommand{\s}{\Sigma}
\newcommand{\con}{  \emph{con}}
\newcommand{\CON}{ \mathrm{CON}}

\newcommand{\ts}{\widetilde{\Sigma}}

\newcommand{\hc}{\circ}
\newcommand{\n}{\{}
\newcommand{\nn}{\}}
\newcommand{\ca}{\eta}

\newcommand{\lt}{\varepsilon}

\newcommand{\fl}{\mathrm{d}}

\newcommand{\cx}{\delta}

\newcommand{\af}{\alpha}

\newcommand{\x}{\langle}
\newcommand{\xx}{\rangle}

\newcommand{\xyd}{\subseteq}
\newcommand{\ty}{ \equiv}
\newcommand{\q}{ \Delta}

\newcommand{\m}{\;\mathrm{mod}\;}

\newtheorem{thm0}{Theorem}

\newtheorem{definition}{Definition}[section]

\newtheorem{lem}{Lemma}[section]
\newtheorem{pro}[lem]{Proposition}

\newenvironment{proof}{\noindent{\bf Proof.}}
{\noindent \ \hfill$\Box$\par}

\begin{document}
\title{\bf{ On the generalized Jacobi Identity of Toda brackets}
}

\author{ Juxin Yang\thanks{Yanqi Lake Beijing Institute of Mathematical Sciences and Applications\,(BIMSA),\;Beijing, 101408,  P.R. China;\;and Shinshu University, 3-1-1 Asahi, Matsumoto, Nagano 390-8621, Japan. EMAIL: yangjuxin@bimsa.com}}

\date{}
\maketitle
\noindent
\begin{abstract}
 In this paper,  we give the  Jacobi identity of Toda brackets indexed by $n,\; (n\geq0)$, this generalizes Toda's Jacobi identity of Toda brackets indexed by 0 shown in \cite{Toda}.

\end{abstract}
\section{Introduction}

 \indent\quad\,Toda bracket is an art of constructing  homotopy liftings   and  homotopy extensions of maps,  it plays a fundamental role on dealing with  composition relations of homotopy classes. We call the formulas of Toda brackets shown in \cite[Proposition\,1.5]{Toda},\;\cite[Formula 3.7)]{Toda} and  \cite[Formula 3.9)\;ii)]{Toda}  Jacobi identities. The Jacobi identities of Toda brackets are critical tools to detect which elements a Toda bracket can contain.\,Toda brackets indexed by $n$ ($n\geq1$) are of great importance to study the unstable phenomena,  for example, to use Toda bracket methods to  consider the image of the boundary homomorphism $\pi_{i+1}(X)\stackrel{\q}\longrightarrow \pi_{i}(F)$ where $F\rightarrow A\rightarrow X$ is a homotopy fibration, it is usually necessary to use Toda brackets indexed by $n$ ($n\geq1$), because the property  $\q\n \af, E^{n}\beta,E^{n}\gamma\nn_{n}\xyd\pm\n \q\af, E^{n-1}\beta,E^{n-1}\gamma\nn_{n-1}$ holds for $n\geq1$ but not for $n=0$; and  Toda brackets indexed by $n$ ($n\geq1$) are also powerful tools  to desuspend homotopy classes. In this paper,  we give the  Jacobi identity  of Toda brackets indexed by $n,\; (n\geq0)$, which generalizes \cite[Prop.\,1.5]{Toda}, the  Jacobi identity  of Toda brackets indexed by $0$. \\\indent
 Here is our main theorem.

\begin{thm0}\label{hao1}

  For the following sequence of spaces and homotopy classes, where  $\af_{1}\hc\s^{n}\af_{2}=\af_{2}\af_{3}=\af_{3}\af_{4}=\af_{4}\af_{5}=0,$\;$n\geq0,$ and  $$\af_{1}\hc \s^{n}\n\af_{2},\af_{3},\af_{4}\nn\ni0, \quad\n\af_{2},\af_{3},\af_{4}\nn\hc\s\af_{5}\ni 0.$$
\[\xymatrix@!C=0.93cm{X&\s^{n}  Y\; & \s^{n}Z &\s^{n}K&\s^{n}L &\s^{n}U \ar"1,3";"1,2" _{\;\s^{n}\af_{2}} \ar"1,4";"1,3" _{
\; \s^{n}\af_{3}}
 \ar"1,5";"1,4" _{
 \;\s^{n}\af_{4}} \ar"1,6";"1,5" _{
\s^{n}\af_{5}}
\ar"1,2";"1,1" _{\af_{1}}
     }\]

\noindent There exist $\xi_{1}\in \n \af_{1},\s^{n}\af_{2},\s^{n}\af_{3}\nn$,
$\xi_{2}\in \n \af_{2},\af_{3},\af_{4}\nn$ and $\xi_{3}\in \n \af_{3},\af_{4},\af_{5}\nn$, such that $$\n\xi_{1},\s^{n+1}\af_{4},\s^{n+1}\af_{5} \nn_{n} +(-1)^{n}\n \af_{1}, \s^{n}\xi_{2},\s^{n+1}\af_{5} \nn_{n}+\n\af_{1},\s^{n}\af_{2},\s^{n}\xi_{3} \nn_{n}\ni0.$$
It is convenient to remember the result  as the following,
\begin{eqnarray}
\notag&
\n\n \af_{1},\s^{n}\af_{2},\s^{n}\af_{3}\nn,\;\s^{n+1}\af_{4},\s^{n+1}\af_{5} \nn_{n}&    \\\notag&
+&    \\\notag&
(-1)^{n}\n \af_{1}, \s^{n} \n \af_{2},\af_{3},\af_{4}\nn,\;\s^{n+1}\af_{5} \nn_{n}&   \\\notag&
+&    \\\notag&
\n\af_{1},\s^{n}\af_{2},\s^{n} \n \af_{3},\af_{4},\af_{5}\nn\nn_{n}&\\\notag&
 |||&  \\\notag&
\;\; 0\,.&
\end{eqnarray}

\end{thm0}

\textbf{Acknowledgement}\;The author  is indebted  to  professor Juno Mukai, Dai Tamaki, Yoshihiro Hirato and Ivy Cao   for   many fruitful conversations on
this project.
And the  author would like to thank professor Dai Tamaki very much for invitation to Shinshu University by his grant JSPS KAKENHI\,(No.\,20K03579).
\section{Preliminaries}

\subsection{Notations}

\indent \quad
 In this paper, all spaces, maps and  homotopy classes
are pointed unless otherwise
stated,  we use $*$ to denote  basepoints and constant maps, $I=[0,1]$ has basepoint 1, and we denote the homotopy classes of constant maps  by 0;  by abuse of notation, following \cite{Toda}, for  homotopy classes $\af_{1}$, $\af_{2},\cdots,\af_{n}$ of constant maps,  we denote ``\;$\af_{1}=0,\;\af_{2}=0,\cdots,\af_{n}=0$\;'' by ``\;$\af_{1}=\af_{2}=\cdots=\af_{n}=0$\;'', although the domains or codomains of them may not be same. Spaces we mean are  pointed spaces satisfying  the inclusions from $\n*\nn$ to them are closed cofibrations, for example, pointed \textit{CW} complexes. We denote the suspension functor by $\s$. There are different ways to define Toda brackets, in this paper, all Toda brackets we mean  follow the original definition in \cite{Toda} given by \textit{Toda}. We use $\mathrm{C}(-)$ to denote the  cone functor, for a map $f:X\rightarrow Y$,  we use $\mathrm{C}^{f}:\mathrm{C}X\rightarrow \mathrm{C}Y$  to denote the  map extended over $f$ by applying  the  cone functor, that is,   $\mathrm{C}^{f}=f\wedge\mathrm{id}_{I}$,  notice $\mathrm{C}X=X\wedge I$ where $I$ has basepoint 1.
For a map $f: X\rightarrow Y$ satisfying $f\simeq*$, we use $\CON(f)$ to denote the set of all of the extensions $\mathrm{C}X\rightarrow Y$ of $f$, that is,
$\CON(f)=\n \;\;\overline{f}: \mathrm{C}X\rightarrow Y\;\, \text{is a map} \;\;|\; \;\overline{f}|_{X}=f \;\nn,$
\noindent a map in $\CON(f)$ is usually denoted by $\con(f)$, essentially, $\con(f)$ is a  homotopy  from $f$ to $*$, if we use the symbol $\con(f)$ directly without special statements, then we always default $\con(f)$ is some fixed map in $\CON(f)$, that is, for convenience, we will not say for some $\con(f)\in\CON(f)$ unless otherwise
stated; and we say a map in $\CON(f)$ a null homotopy of $f$. Let $f:X\wedge Y\rightarrow Z$ be a map, we use the obvious notation $f(-\wedge y_{0})$ to denote the map $X\rightarrow Z$ given by $x\mapsto f(x\wedge y_{0})$, for $f(x_{0} \wedge -)$, similarly. And  similar notations $g(-, y_{0})$, $g(x_{0}, -)$ for a map $g:X\times Y\rightarrow Z$ are also used. For $\s X=X\wedge([0,1]/\n0,1\nn)$, we use $[k]: \s X\rightarrow \s X$ to denote the map of degree $k$ where $k\in\z$, in particular, we use $[(-1)^{n}]$ to  denote the map of degree $(-1)^{n}$, speaking
in detail, $[0]=*$, $[1]=\mathrm{id}$, $[-1]$ is given by $x\wedge t\mapsto x\wedge(1-t)$, further  $[(-1)^{m}\ell]$ is $[(-1)^{m}]+[(-1)^{m}]+\cdots+[(-1)^{m}]$, ($\ell$ times, $\ell\in\z_{+}$, $m\in\n0,1\nn$),
 by abuse of notation, the homotopy class of $[k]$ is still denoted by $[k]$. We use $\tau(X,Y): X\wedge Y\rightarrow Y\wedge X$ to denote the exchanging coordinates map given by $\tau(X,Y)(x\wedge y)= y\wedge x$.

\subsection{Some fundamental knowledge}
\subsubsection{The Suspension Tilde Operator}
\indent \quad The main use of  the suspension tilde operator $\ts^{n}$ is  to handle  the classical homemorphism $\s^{n}\mathrm{C}X\rightarrow \mathrm{C}\s^{n}X$, see Lem.\,\ref{tsyl}\,(2) and Prop.\,\ref{djdy}.

\begin{definition}(\cite{dd}) Let $X,Y$ be spaces and $n\geq1$, for a map $F:\mathrm{C}X\rightarrow Y$, define $\ts^{n}F$ to be
   $\ts^{n}F:\mathrm{C}\s^{n}X\rightarrow \s^{n}Y$, that is,
    $$\ts^{n}F:X\wedge S^{n}\wedge I\rightarrow Y\wedge S^{n},$$ it is given by $$x\wedge a\wedge t\mapsto F(x\wedge t)\wedge a,\quad \;\;\text{for any\;}x\wedge a\wedge t\in X\wedge S^{n}\wedge I.$$
  And  define $\ts^{0}F=F.$
\end{definition}

\indent For the following lemma, (1) and (2)   are from \cite[p.\,17]{dd}, and  (3) is immediately got by the definition of the operator $\ts^{n}$.\\
\indent We regard $\s^{0}=\mathrm{id}$ and $\mathrm{id}_{X}\wedge \tau(S^{0}, I) =\mathrm{id}_{\mathrm{C}X}$.
\begin{lem}\label{tsyl} For a map $F:\mathrm{C}X\rightarrow Y$,
    \begin{itemize}
        \item[\rm(1)]  $\ts^{m}(\ts^{n}F)=\ts^{m+n}F$,\; for any non-negative $m,n$;
         \item[\rm(2)]  $\ts^{n}F=\s^{n}F\hc (\mathrm{id}_{X}\wedge \tau(S^{n}, I))$,\; $n\geq0$;
          \item[\rm(3)]   if $g:X'\rightarrow X$ is a map, then $$\ts^{n}(F\hc \mathrm{C}^{g})=(\ts^{n}F)\hc \mathrm{C}^{\s^{n}g}, \;n\geq0;$$ if $G:\mathrm{C}L\rightarrow X$ is a map, then
          $$ (\ts^{n}(F\hc \mathrm{C}^{G}))\hc (\mathrm{id}_{L}\wedge\mathrm{id}_{I}\wedge \tau(S^{n}, I))=(\ts^{n}F)\hc \mathrm{C}^{\ts^{n}G},\;n\geq0,$$
 $$ (\s^{n}(F\hc \mathrm{C}^{G}))\hc (\mathrm{id}_{L}\wedge \tau(S^{n}, I\wedge I))=(\ts^{n}F)\hc \mathrm{C}^{\ts^{n}G},\;n\geq0.$$

    \end{itemize}
\end{lem}
\subsubsection{The Separation Element}
\indent\quad  To preserve  the order of the maps when we are dealing with Toda brackets,   we denote the notation $d(A,B)$  in \cite[p.\,19]{dd} by $\fl(B,A)$, denote the notation $\cx(A,B)$  in \cite{dd} by $\cx(B,A)$. Hereafter, $S^{1}=[0,1]/\n0,1\nn$, we denote $q(t)$ by $t$ for simplicity where $q:[0,1]\rightarrow S^{1}$ is the projection.

\begin{definition}(\cite{dd}) Suppose $X,Y$ are spaces, and $f,g:\mathrm{C}X\rightarrow Y$ are two maps, $f\,|\,_{X}=g\,|\,_{X}$.  The separation element of $(f,g)$ which is denoted by $\fl(f,g)$ is defined to be $\fl(f,g): \s X=X\wedge ([0,1]/\n\,0,1\,\nn)  \rightarrow Y$  which is given by  $$x\wedge t \mapsto g(x \wedge (1-2t)),\quad t\in[0,1/2],$$ $$x\wedge t \mapsto f(x \wedge (2t-1)),\quad t\in[1/2,1], $$
and define $\cx(f,g)$ to be the homotopy class of $\fl(f,g)$.
\end{definition}

Sometimes, in order to avoid the inconvenience caused by too many symbols $``(\,\;, \;)$'', we also denote $\fl(f,g)$ by $\fl\n f,g\nn$. For example, we write $$\fl(f_{1}\hc \con (f_{2}f_{3})\hc \mathrm{C}^{f_{4}},\;\con (f_{5}f_{6})\hc \mathrm{C}^{f_{7}}\hc \s(f_{8}f_{9})) $$ by
$$\fl\n f_{1}\hc \con (f_{2}f_{3})\hc \mathrm{C}^{f_{4}},\;\con (f_{5}f_{6})\hc \mathrm{C}^{f_{7}}\hc \s(f_{8}f_{9}) \nn.$$

\indent The following lemma is from \cite[Lem.\,2.8]{dd} except for the second part of (3) which is easy to check.
\begin{lem}\label{flyl}Suppose $f,g,h:\mathrm{C}X\rightarrow Y$ are maps, $f\,|\,_{X}=g\,|\,_{X}=h\,|\,_{X}$.
\begin{itemize}

    \item[\rm(1)] $\fl(f,g)=-\fl(g,f)$.
    \item[\rm(2)] $\cx(f,g)+\cx(h,f)=\cx(g,h)$.
    \item[\rm(3)] $\cx(f,g)=0$ if and only if $f\simeq g\;\mathrm{rel}\;X$. \\ Successively, if $f\simeq f'\;\mathrm{rel}\;X$ and $g\simeq g'\;\mathrm{rel}\;X$ , then, $\fl(f,g)\simeq \fl(f',g') .$
    \item[\rm(4)] For maps $r: W\rightarrow X$ and $s: Y\rightarrow Z$, the following relations hold,$$ \fl(f,g)\hc\s r= \fl(f\hc \mathrm{C}^{r},g\hc \mathrm{C}^{r}),\quad s\hc\fl(f,g)= \fl(s\hc f, s\hc g).$$
     \item[\rm(5)]  $\s^{n}\fl(f,g)\simeq \fl(\ts^{n}f,\ts^{n} g)\hc[(-1)^{n} ]$,\;$n\geq0$.
\end{itemize}

\end{lem}

For a map $f:\mathrm{C}X\rightarrow Y$, let $f^{+\frac{1}{2}}: \mathrm{C}X\rightarrow Y$ be the map given by  \[f^{+\frac{1}{2}}(x\wedge t)=\begin{cases}
f(x\wedge (2t-1)),\;&\text{    $t\in[\frac{1}{2},1],$}\\
f(x\wedge 0),\;&\text{$t\in[0,\frac{1}{2}].$}
\end{cases}\]
\noindent By use of this notation, we give the following lemma.

\begin{lem}\label{efzyyl} Let $f,g:\mathrm{C}X\rightarrow Y$
    be maps, then $$\fl(f,g)\simeq \fl(f^{+\frac{1}{2}},g)\simeq \fl(f,g^{+\frac{1}{2}})\simeq \fl(f^{+\frac{1}{2}},g^{+\frac{1}{2}}).$$
\end{lem}
\begin{proof}
We only need to show $f^{+\frac{1}{2}}\simeq f\;\mathrm{rel}\;\,X$, then this lemma holds by Lem.\,\ref{flyl}\,(3). \\
\indent As a matter of fact, we have a function  $H: (X\wedge I)\times I\rightarrow Y$ which is given by
\[H(x\wedge t, r)=\begin{cases}
f(\;x\wedge \frac{2}{r+1}(t-\frac{1-r}{2}) \;),\;&\text{    $t\in[\frac{1-r}{2},1],$}\\
f(x\wedge 0),\;&\text{$t\in[0,\frac{1-r}{2}].$}
\end{cases}\]
\noindent Therefore, $H(-,0)=f^{+\frac{1}{2}}$,  $H(-,1)=f$, $H(*\wedge 0,-)=*$. We notice $$\n (x\wedge \frac{1-r}{2},r)\;\;|\;\; x\in X, r\in[0,1]\nn$$
is a closed set in $(X\wedge I)\times I$, so, $H$ is  continuous, that is, $H$ is a map. Thus, $H$ is  a homotopy from $f^{+\frac{1}{2}}$ to $f$. To check $H$ is  a homotopy relative to $X$, let $t=0$,  we have $H(x\wedge 0, r)=f(x\wedge 0)$ for all $r$. So $H$ is truly a homotopy relative to $X.$ Hence our lemma holds.

\end{proof}

\;\\\indent
Recall for a map $f: X\rightarrow Y$ satisfying $f\simeq*$, we use $\CON(f)$ to denote the set of all of the extensions $\mathrm{C}X\rightarrow Y$ of $f$, that is,
$$\CON(f)=\n \;\;\overline{f}: \mathrm{C}X\rightarrow Y\;\, \text{is a map} \;\;\;\;|\;\;\; \;\overline{f}|_{X}=f \;\nn,$$
\noindent a map in $\CON(f)$ is usually denoted by $\con(f)$.\\\indent
The following proposition is due to \cite[Prop.\,2.11]{dd}, \cite[Lem.\,2.4]{dd} and \cite[p.\,22]{dd}.
\begin{pro}(\cite{dd})\label{djdy}
    For the following sequence of spaces and homotopy classes, where  $\af_{1}\hc \s^{n}\af_{2}=\af_{2}\af_{3}=0,$
\[\xymatrix@!C=0.83cm{ &W&  \s^{n}X\; & \s^{n}Y  &\s^{n}Z,&
       \ar"1,4";"1,3" _{\;\s^{n}\af_{2}} \ar"1,5";"1,4" _{
 \s^{n}\af_{3}}
\ar"1,3";"1,2" _{\af_{1}}
     }\]
 the  Toda bracket  $\n\af_{1},\s^{n}\af_{2},\s^{n}\af_{3}\nn_{n}$  is equal to the collection of  the homotopy classes of form $$ \cx (\,f_{1}\hc\ts^{n}\con(f_{2}f_{3})\;,\; \con(f_{1}\hc\s^{n}f_{2})\hc \mathrm{C}^{\s^{n}f_{3}}\,), $$\noindent
where   $f_{i}$ can be any map in $\af_{i},$ after choosing $f_{i}$, $\con(f_{2}f_{3})$  can be any map in $\CON(f_{2}f_{3})$, and  $\con(f_{1}\hc\s^{n}f_{2})$  can be any map in $\CON(f_{1}\hc \s^{n}f_{2})$.
\end{pro}
\section{The Jacobi Identity of Toda Brackets Indexed by $n$}
\begin{lem}\label{ndyl}
     For the following sequence of spaces and homotopy classes where  $\af_{1}\af_{2}=\af_{3}\af_{4}=0,$
\[\xymatrix@!C=0.53cm{X&  Y\; & Z &K&L, \ar"1,3";"1,2" _{\af_{2}} \ar"1,4";"1,3" _{
 \af_{3}}
 \ar"1,5";"1,4" _{
 \af_{4}}
\ar"1,2";"1,1" _{\af_{1}}
     }\]
then the relation
$$\con(f_{1}f_{2})\hc \mathrm{C}^{f_{3}f_{4}}\simeq f_{1}f_{2}\hc \con (f_{3}f_{4})\;\;\mathrm{rel}\;L$$
holds  for any $f_{k}\in\af_{k}$,  $\con(f_{1}f_{2})\in \CON(f_{1}f_{2})$,  $\con(f_{3}f_{4})\in \CON(f_{3}f_{4})$.
\end{lem}

\begin{proof}
For any given $f_{k}\in\af_{k}$,  $\con(f_{1}f_{2})\in \CON(f_{1}f_{2})$,  $\con(f_{3}f_{4})\in \CON(f_{3}f_{4})$,
\noindent we construct a homotopy $F: (\mathrm{C}L)\times I\rightarrow X$\; $\mathrm{rel}\;L$.
Define $$F(\ell \wedge s, t)=\con(f_{1}f_{2})\hc \mathrm{C}^{(\con(f_{3}f_{4}))(-\wedge st)}(\ell\wedge s(1-t))$$
where $\ell\wedge s\in L\wedge I=\mathrm{C}L,\;t\in I.$ Then,  $F(* \wedge 0,-)=*$, and $$F(\ell \wedge s, 0)=\con(f_{1}f_{2})\hc \mathrm{C}^{f_{3}f_{4}},$$
   $$F(\ell \wedge s,1)=f_{1}f_{2}\hc \con (f_{3}f_{4}),$$  \noindent Hence $F$ is truly a homotopy. To check this homotopy is relative to $L$, let $s=0$, we have
   $$F(\ell \wedge 0, 0)=F(\ell \wedge 0, t)\;\,\text{ for any}\; t.$$\noindent  Hence the result holds.
\end{proof}

\begin{lem}\label{qryl} Suppose $X$ and $Y$ are spaces, $f,g:X\rightarrow Y$ are maps, $f\simeq g\simeq*$. Then, there exist $Q\in \CON (f)$, $R\in \CON (g)$, such that $$R(x\wedge t)=Q(x\wedge\, (1+t)/2 ),\;\,\text{for all}\;\,x\wedge t\in X\wedge I=\mathrm{C}X.$$

\end{lem}
\begin{proof}
  By the given, there exist homotopies $M,N:X\times I\rightarrow Y$, such that $$
  M(-,0)=f, \quad M(-,1)=N(-,0)=g,\quad N(-,1)=*,$$ $$ M(*,-)=N(*,-)=*.$$

\noindent Then, we have  a homotopy $P:X\times I\rightarrow Y$ from $f$ to $*$ which is given by$$ (x,t)\mapsto M(x,2t),\; t\in[0,1/2],$$$$(x,t)\mapsto N(x,2t-1),\; t\in[1/2,1]. $$
  Successively, we obtain a map $Q:\mathrm{C}X\rightarrow Y$ given by $ x\wedge t\mapsto P(x,t)$ satisfying $Q(-\wedge 0)=f.$ Then, we get a map $R:\mathrm{C}X\rightarrow Y$ given by $x\wedge t\mapsto Q(x\wedge\, (1+t)/2)$, it satisfies $$R(-\wedge0)= Q(-\wedge\,1/2)=P(-,1/2)=M(-,1)=g.$$Hence the result holds.
\end{proof}
\begin{lem}\label{zb}
     For the following sequence of spaces and homotopy classes, where  $\af_{1}\hc\s^{n}\af_{2}=\af_{2}\af_{3}=\af_{3}\af_{4}=\af_{4}\af_{5}=0,$  $n\geq0,$ and  $$\af_{1}\hc \s^{n}\n\af_{2},\af_{3},\af_{4}\nn\ni0, \quad\n\af_{2},\af_{3},\af_{4}\nn\hc\s\af_{5}\ni 0.$$
\[\xymatrix@!C=0.93cm{X&\s^{n}  Y\; & \s^{n}Z &\s^{n}K&\s^{n}L &\s^{n}U \ar"1,3";"1,2" _{\;\s^{n}\af_{2}} \ar"1,4";"1,3" _{
\; \s^{n}\af_{3}}
 \ar"1,5";"1,4" _{
 \;\s^{n}\af_{4}} \ar"1,6";"1,5" _{
\s^{n}\af_{5}}
\ar"1,2";"1,1" _{\af_{1}}
     }\]

\noindent There exists $\xi_{2}\in \n \af_{2},\af_{3},\af_{4}\nn$, such that $\af_{1}\hc\s^{n} \xi_{2}=\xi_{2}\hc\s\af_{5}=0$.
\end{lem}
\begin{proof}
   By the given, there exists $\lt\in\n\af_{2},\af_{3},\af_{4}\nn$, such that $\af_{1}\hc \s^{n}\lt=0$.
Since $\n\af_{2},\af_{3},\af_{4}\nn\hc\s\af_{5}\ni 0$ and $\af_{4}\af_{5}=0$,
then, $$\lt\hc\s\af_{5}\ty0\m \af_{2}\hc [\s L,Z]\hc\s\af_{5}.$$
   So, there exists $\ca\in[\s L, Z]$, such that $\lt\hc\s\af_{5}=\af_{2}\ca\hc\s\af_{5}$. Let $$\xi_{2}=\lt+\af_{2}\ca\hc[-1],$$
\noindent then $\xi_{2}\hc\s\af_{5}=0$, and $\af_{1}\hc\s^{n}\xi_{2}=0.$
\end{proof}

\begin{pro}\label{hu}
     For the following sequence of spaces and homotopy classes, where  $\af_{1}\hc\s^{n}\af_{2}=\af_{2}\af_{3}=\af_{3}\af_{4}=\af_{4}\af_{5}=0,$\;$n\geq0,$ and  $$\af_{1}\hc \s^{n}\n\af_{2},\af_{3},\af_{4}\nn\ni0, \quad\n\af_{2},\af_{3},\af_{4}\nn\hc\s\af_{5}\ni 0.$$
\[\xymatrix@!C=0.93cm{X&\s^{n}  Y\; & \s^{n}Z &\s^{n}K&\s^{n}L &\s^{n}U \ar"1,3";"1,2" _{\;\s^{n}\af_{2}} \ar"1,4";"1,3" _{
\; \s^{n}\af_{3}}
 \ar"1,5";"1,4" _{
 \;\s^{n}\af_{4}} \ar"1,6";"1,5" _{
\s^{n}\af_{5}}
\ar"1,2";"1,1" _{\af_{1}}
     }\]

\noindent There exist $\xi_{1}\in \n \af_{1},\s^{n}\af_{2},\s^{n}\af_{3}\nn$,
$\xi_{2}\in \n \af_{2},\af_{3},\af_{4}\nn$ and $\xi_{3}\in \n \af_{3},\af_{4},\af_{5}\nn$, such that $$\n\xi_{1},\s^{n+1}\af_{4},\s^{n+1}\af_{5} \nn_{n} +(-1)^{n}\n \af_{1}, \s^{n}\xi_{2},\s^{n+1}\af_{5} \nn_{n}+\n\af_{1},\s^{n}\af_{2},\s^{n}\xi_{3} \nn_{n}\ni0.$$
It is convenient to remember the result  as the following,
\begin{eqnarray}
\notag&
\n\n \af_{1},\s^{n}\af_{2},\s^{n}\af_{3}\nn,\;\s^{n+1}\af_{4},\s^{n+1}\af_{5} \nn_{n}&    \\\notag&
+&    \\\notag&
(-1)^{n}\n \af_{1}, \s^{n} \n \af_{2},\af_{3},\af_{4}\nn,\;\s^{n+1}\af_{5} \nn_{n}&   \\\notag&
+&    \\\notag&
\n\af_{1},\s^{n}\af_{2},\s^{n} \n \af_{3},\af_{4},\af_{5}\nn\nn_{n}&\\\notag&
 |||&  \\\notag&
\;\; 0\,.&
\end{eqnarray}

\end{pro}
\begin{proof} This proof is divided into 4 steps. The \textit{2nd, 3rd} and \textit{4th} steps parallel to the second half of the proof of \cite[Prop. 4.2 ii)]{T59}, but our situation is more complicated, and the explanation of us gives  more details. \\
\indent\textbf{STEP\,1.}\; We show there are maps $p_{1},p_{2},p_{3}$, their homotopy classes are the desired $\xi_{1},\xi_{2},\xi_{3}$ respectively, and show some properties of them.\\ \indent
 By Lem.\,\ref{zb}, we know   there  exists $$\xi_{2}\in \n \af_{2},\af_{3},\af_{4}\nn,\;\text{such that}\; \af_{1}\hc\s^{n} \xi_{2}=\xi_{2}\hc\s\af_{5}=0.$$
By Prop.\,\ref{djdy}, we know
there  exist $f_{k}\in\af_{k}\, (k=2,3,4)$ and $\con(f_{3}f_{4}), \con(f_{3}f_{3})$,  such that $$p_{2}:=\fl\n f_{2}\hc \con(f_{3}f_{4}),\;\con(f_{2}f_{3})\hc \mathrm{C}^{f_{4}}\nn\in \xi_{2}.$$
 We take $f_{1}\in\af_{1}$, by Lem.\,\ref{flyl}\,(1),(4),(5) and Lem.\,\ref{tsyl}\,(2),(3), we have
\begin{eqnarray}
&
&\notag f_{1}\hc\s^{n} p_{2}  \\\notag&
\simeq& f_{1}\hc \fl\n\ts^{n}(f_{2}\con(f_{3}f_{4})),\;\ts^{n}(\con(f_{2}f_{3})\hc \mathrm{C}^{f_{4}})\nn\hc [(-1)^{n}]  \\\notag&
=& \fl\n f_{1}\hc\ts^{n}(f_{2}\con(f_{3}f_{4})),\;f_{1}\hc\ts^{n}(\con(f_{2}f_{3})\hc \mathrm{C}^{f_{4}})\nn\hc[(-1)^{n}]  \\\notag&
=&\fl\n f_{1}\hc\s^{n}f_{2}\hc\s^{n}\con(f_{3}f_{4})\hc\tau_{1},\;f_{1}\hc\ts^{n}\con(f_{2}f_{3})\hc \mathrm{C}^{\s^{n}f_{4}}\nn\hc[(-1)^{n}]  \\\notag&
=& \fl\n f_{1}\hc\s^{n}f_{2}\hc\con(\s^{n}f_{3}\hc\s^{n}f_{4}),f_{1}\hc\ts^{n}\con(f_{2}f_{3})\hc \mathrm{C}^{\s^{n}f_{4}}\nn\hc[(-1)^{n}]  \\\notag&
=&  \fl\n f_{1}\hc\ts^{n}\con(f_{2}f_{3})\hc \mathrm{C}^{\s^{n}f_{4}},\;f_{1}\hc\s^{n}f_{2}\hc\con(\s^{n}f_{3}\hc\s^{n}f_{4})\nn\hc[(-1)^{n+1}],
\end{eqnarray}

\noindent here $\tau_{1}=\mathrm{id}_{L}\wedge\tau(I,S^{n})$, and we notice $\con(\s^{n}f_{3}\hc\s^{n}f_{4})$ can be taken as $\s^{n}\con(f_{3}f_{4})\hc\tau_{1}$.
\noindent By Lem.\,\ref{ndyl}, we know $$f_{1}\hc\s^{n}f_{2}\hc\con(\s^{n}f_{3}\hc\s^{n}f_{4})\simeq \con(f_{1}\hc\s^{n}f_{2})\hc \mathrm{C}^{\s^{n}f_{3}}\hc \mathrm{C}^{\s^{n}f_{4}}\;\mathrm{rel}\;\s^{n}L.$$
Then, by Lem.\,\ref{flyl}\,(3),(4) we have,
\begin{eqnarray}
&
&\notag \fl\n f_{1}\hc\ts^{n}\con(f_{2}f_{3})\hc \mathrm{C}^{\s^{n}f_{4}},\;f_{1}\hc\s^{n}f_{2}\hc\con(\s^{n}f_{3}\hc\s^{n}f_{4})\nn\hc[(-1)^{n+1}] \\\notag&
\simeq&  \fl\n f_{1}\hc\ts^{n}\con(f_{2}f_{3})\hc \mathrm{C}^{\s^{n}f_{4}},\;\con(f_{1}\hc\s^{n}f_{2})\hc \mathrm{C}^{\s^{n}f_{3}}\hc \mathrm{C}^{\s^{n}f_{4}}\nn\hc[(-1)^{n+1}]   \\\notag&
=& \fl\n f_{1}\hc\ts^{n}\con(f_{2}f_{3}),\;\con(f_{1}\hc\s^{n}f_{2})\hc \mathrm{C}^{\s^{n}f_{3}}\nn\hc \s^{n+1}f_{4}\hc[(-1)^{n+1}].
\end{eqnarray}
So,
$$f_{1}\hc\s^{n}p_{2}\simeq\fl\n f_{1}\hc\ts^{n}\con(f_{2}f_{3}),\;\con(f_{1}\hc\s^{n}f_{2})\hc \mathrm{C}^{\s^{n}f_{3}}\nn\hc \s^{n+1}f_{4}\hc[(-1)^{n+1}].$$
\noindent By Prop.\,\ref{djdy} we know there exists $\xi_{1}\in\n\af_{1},\s^{n}\af_{2},\s^{n}\af_{3}\nn_{n}$
such that\\\\ \centerline{$\xi_{1}
\ni\fl\n f_{1}\hc\ts^{n}\con(f_{2}f_{3}),\,\con(f_{1}\hc\s^{n}f_{2})\hc \mathrm{C}^{\s^{n}f_{3}}\nn:=p_{1},$}
then,\begin{subequations}
\begin{align}
p_{1}\hc \s^{n+1}f_{4}\simeq  f_{1}\hc \s^{n}p_{2}\hc[(-1)^{n+1}]\simeq*, \label{yfh1}
\end{align}
\end{subequations}
(recall $f_{1}\hc\s^{n}p_{2}\in\af_{1}\hc\s^{n}\xi_{2}=0$).
 Similarly, for some $f_{5}\in\af_{5}$, there exists $\xi_{3}\in\n\af_{3},\af_{4},\af_{5}\nn$ such that $\xi_{3}\ni \fl\n f_{3}\hc\con(f_{4}f_{5}),\; \con(f_{3}f_{4}) \hc \mathrm{C}^{f_{5}}   \nn:=p_{3},\;\text{and}$
\begin{subequations}
\begin{align}
f_{2}p_{3}\simeq p_{2}\hc\s f_{5}\hc[-1]\simeq*. \label{yfh2}
\end{align}
\end{subequations}
\indent\textbf{STEP\,2.}\; We show some properties of the null homotopies used to construct the  desired representative homotopy classes. \\\\
\indent  For a map $f$, let $dom(f)$  be the domain of $f$. We take
  the null homotopies of the  constant maps as the following, $$\mathscr{A}=\con(f_{1}\s^{n}f_{2}),\; dom(\mathscr{A})=\mathrm{C}\s^{n}Z=Z\wedge S^{n}\wedge I,$$
$$\mathscr{B}=\con(f_{2}f_{3}),\; dom(\mathscr{B})=\mathrm{C}K=K\wedge  I,$$
$$\mathscr{C}=\con(f_{3}f_{4}),\;dom(\mathscr{C})=\mathrm{C}L=L\wedge I ,$$
$$\mathscr{D}=\con(f_{4}f_{5}),\;dom(\mathscr{D})=\mathrm{C}U=U\wedge I.$$
And let
$\mathscr{A}_{t}=\mathscr{A}(-\wedge t)$,
$\mathscr{B}_{t}=\mathscr{B}(-\wedge t)$, $\mathscr{C}_{t}=\mathscr{C}(-\wedge t)$,
$\mathscr{D}_{t}=\mathscr{D}(-\wedge t)$,  where $t\in I$ is the last coordinate corresponding the last smash factors  shown in the above domains.
\noindent By Lem.\,\ref{qryl}, there exist maps $$\con(p_{1}\hc\s^{n+1} f_{4})\;\text{ and}\;\; \con(f_{1}\hc\s^{n}((-1)^{n+1}p_{2}))$$ both with domain $\mathrm{C}\s^{n+1}L=L\wedge S^{n+1}\wedge I$, and they satisfy
\begin{subequations}
\begin{align}
\con(f_{1}\hc\s^{n}(-1)^{n+1}p_{2}))(x\wedge t)=\con(p_{1}\hc\s^{n+1} f_{4})(x\wedge (1+t)/2)\;\label{cltl1}
\end{align}
\end{subequations}

\noindent$\text{for all}\; x\in L\wedge S^{n+1}, t\in I.$ Similarly,
there exist maps $\con(f_{2}\hc(-1)^{n}p_{3})$ and $\con((-1)^{n+1}p_{2}\hc (-1)^{n}\s f_{5})$ both with domain $\mathrm{C}\s U=U\wedge S^{1}\wedge I$ such that
\begin{subequations}
\begin{align}
\con((-1)^{n+1}p_{2}\hc (-1)^{n}\s f_{5})(x\wedge t)=\con(f_{2}\hc(-1)^{n}p_{3})(x\wedge\, (1+t)/2)\;\label{cltl2}
\end{align}
\end{subequations}
\noindent
$\text{for all}\; x\in U\wedge S^{1}, t\in I.\;$\noindent
We choose the maps $$\mathcal{G,\; G'}:\mathrm{C}\s^{n+1}L=L\wedge S^{n+1}\wedge I\rightarrow X,$$
$$\mathscr{F,\; F'}:\mathrm{C}\s^{n+1}U=U\wedge S^{n+1}\wedge I\rightarrow \s^{n}Y,$$\noindent
as the following,
 $$\mathcal{G}=\con(f_{1}\hc\s^{n}(-1)^{n+1}p_{2}),\quad\mathcal{G}'=\con(p_{1}\hc\s^{n+1}f_{4}),$$
$$\mathscr{F}=\ts^{n}\con((-1)^{n+1}p_{2}\hc (-1)^{n}\s f_{5}),\quad\mathscr{F}'=\ts^{n}\con(f_{2}\hc(-1)^{n}p_{3}),$$
\noindent and let
$\mathscr{F}_{t}=\mathscr{F}(-\wedge t)$,
$\mathscr{F}'_{t}=\mathscr{F}'(-\wedge t)$, $\mathcal{G}_{t}=\mathcal{G}(-\wedge t)$,
$\mathcal{G}'_{t}=\mathcal{G}'(-\wedge t)$,  where $t\in I$ is the last coordinate corresponding the last smash factors  shown in the above domains. Then, by formula\,\ref{cltl1} and formula\,\ref{cltl2}, we know
\begin{subequations}
\begin{align}
\mathcal{G}'_{(1+t)/2}=\mathcal{G}_{t}\;,\quad\mathscr{F}'_{(1+t)/2}=\mathscr{F}_{t}\;. \label{gx}
\end{align}
\end{subequations}
\indent\textbf{STEP\,3.}
We construct maps $\psi_{k}$, such that the homotopy classes $[\psi_{k}]$ can be contained in $ \Psi_{k},\; (k=1,2,3)$ where
$$\Psi_{1}=\n\xi_{1},\;\s^{n+1}\af_{4},\;\s^{n}((-1)^{n}\s\af_{5})\nn_{n},\; $$
$$\Psi_{2}=-\n \af_{1},\; \s^{n}((-1)^{n+1}\xi_{2})\,,\;\s^{n}((-1)^{n}\s\af_{5})\nn_{n},$$$$\Psi_{3}=\n\af_{1},\s^{n}\af_{2},\s^{n}((-1)^{n}\xi_{3}) \nn_{n}.$$
\indent We take $\con((-1)^{n+1}\s f_{4}\hc\s f_{5})=[(-1)^{n+1}]\hc\ts\con( f_{4}f_{5})=[(-1)^{n+1}]\hc\ts \mathscr{D}$,  by calculation, we can obtain directly the following,$$\ts^{n}\con((-1)^{n+1}\s f_{4}\hc\s f_{5})= \tau_{2}\hc\ts^{n+1} \mathscr{D} $$
\noindent
where $\tau_{2}=
\mathrm{id}_{K}\wedge ((-1)^{n+1}\mathrm{id}_{S^{1}})\wedge \mathrm{id}_{S^{n}}.$
\noindent
Let $\tau_{3}=\mathrm{id}_{K}\wedge \tau(S^{1},S^{n})$, then,
by Prop.\,\ref{djdy}, we can take
$$\psi_{1}=\fl( p_{1}\tau_{3}\tau_{2}\hc\ts^{n+1}\mathscr{D}, \;\mathcal{G}'\hc \mathrm{C}^{\s^{n+1}f_{5}} ),\;\text{such that}\; [\psi_{1}]\in\Psi_{1}.$$
Successively,
$$ p_{1}\tau_{3}\tau_{2}\hc\ts^{n+1}\mathscr{D}=\fl\n f_{1}\hc\ts^{n}\con(f_{2}f_{3}),\,\con(f_{1}\s^{n}f_{2})\hc \mathrm{C}^{\s^{n}f_{3}}\nn\hc\tau_{3}\tau_{2}\hc\ts^{n+1}\con( f_{4}f_{5}),$$
\noindent
further we have, for $u\wedge s\wedge a\wedge t\in U\wedge S^{1}\wedge S^{n}\wedge I,\;(S^{1}=[0,1]/\n0,1\nn ),$\begin{subequations}
\begin{align}(p_{1}\tau_{3}\tau_{2}\hc\ts^{n+1}\mathscr{D})(u\wedge s\wedge a\wedge t)=\begin{cases}
(\;\mathscr{A}_{1-2s}(\s^{n}(f_{3}\mathscr{D}_{t}))\;)(u\wedge a),\;&\text{    $s\in[0,1/2],$}\label{p1d}\\
(\;f_{1}(\s^{n}(\mathscr{B}_{2s-1}\mathscr{D}_{t}))\;)(u\wedge a),\;&\text{$\;\,s\in [1/2,1],$}
\end{cases}
\end{align}
\end{subequations}
\noindent
in the case $n$ is odd, and
\begin{subequations}
\begin{align}(p_{1}\tau_{3}\tau_{2}\hc\ts^{n+1}\mathscr{D})(u\wedge s\wedge a\wedge t)=\begin{cases}(\;f_{1}(\s^{n}(\mathscr{B}_{1-2s}\mathscr{D}_{t}))\;)(u\wedge a)
,\;&\text{    $s\in[0,1/2],$}\label{p1dp}\\
(\;\mathscr{A}_{2s-1}(\s^{n}(f_{3}\mathscr{D}_{t}))\;)(u\wedge a)
,\;&\text{$\;\,s\in [1/2,1],$}
\end{cases}
\end{align}
\end{subequations}
\noindent
in the case $n$ is even.  By Lem.\,\ref{flyl}\,(1),
 we can take
$$\psi_{2}'=\fl( f_{1}\hc\mathscr{F},\; \mathcal{G}\hc \mathrm{C}^{\s^{n}((-1)^{n}\s f_{5})}),\;\text{such that}\; [\psi_{2}']\in-\Psi_{2},$$
\noindent then we can choose a map $\psi_{2}''$ as the following,

$$\psi_{2}''=\fl( \mathcal{G}\hc \mathrm{C}^{\s^{n}((-1)^{n}\s f_{5})}, f_{1}\hc\mathscr{F}),\;\;[\psi_{2}'']=[-\psi_{2}']\in\Psi_{2}. $$
\noindent By Lem.\,\ref{efzyyl}, we can take
$$\psi_{2}=\fl\n\;(\mathcal{G}\hc \mathrm{C}^{\s^{n}((-1)^{n}\s f_{5})})^{+\frac{1}{2}},\;(f_{1}\hc\mathscr{F})^{+\frac{1}{2}}\;\nn,\;\text{such that}\; \psi_{2}\simeq \psi_{2}''\,,$$\noindent
successively,  $[\psi_{2}]\in\Psi_{2}.$\;
Similarly to the above, we can take $$\psi_{3}=\fl\n f_{1}\hc\mathscr{F}',\; \mathscr{A}\hc \mathrm{C}^{\s^{n}((-1)^{n}p_{3})}\nn,\;\;\text{such that}\;\,[\psi_{3}]\in\Psi_{3},$$  and
$$ \mathscr{A}\hc \mathrm{C}^{\s^{n}((-1)^{n}p_{3})}= \mathscr{A}\hc \mathrm{C}^{\s^{n}(-1)^{n}\fl\n f_{3}\hc\con(f_{4}f_{5}),\; \con(f_{3}f_{4}) \hc \mathrm{C}^{f_{5}}   \nn},$$
its domain is $\mathrm{C}\s^{n+1}U=U\wedge S^{1}\wedge S^{n}\wedge I$,
then, by direct calculation, we have  for any $u\wedge s\wedge a\wedge t\in U\wedge S^{1}\wedge S^{n}\wedge I,\;(S^{1}=[0,1]/\n0,1\nn ),$

\begin{subequations}
\begin{align}(\mathscr{A}\hc \mathrm{C}^{\s^{n}((-1)^{n}p_{3})})(u\wedge s\wedge a\wedge t)=\begin{cases}(\;\mathscr{A}_{t}\s^{n}(\mathscr{C}_{1-2s}f_{5})\;)(u\wedge a)
,\;&\text{    $s\in[0,1/2]$}\label{ap3}\\
 (\;\mathscr{A}_{t}(\s^{n}(f_{3}\mathscr{D}_{2s-1}))\;)(u\wedge a),\;&\text{$\;\,s\in [1/2,1],$}
\end{cases}
\end{align}
\end{subequations}
\noindent in the case $n$ is odd, and
\begin{subequations}
\begin{align}(\mathscr{A}\hc \mathrm{C}^{\s^{n}((-1)^{n}p_{3})})(u\wedge s\wedge a\wedge t)=\begin{cases}
(\;\mathscr{A}_{t}(\s^{n}(f_{3}\mathscr{D}_{1-2s}))\;)(u\wedge a),\;&\text{    $s\in[0,1/2]$}\label{ap3p}\\
(\;\mathscr{A}_{t}\s^{n}(\mathscr{C}_{2s-1}f_{5})\;)(u\wedge a) ,\;&\text{$\;\,s\in [1/2,1],$}
\end{cases}
\end{align}
\end{subequations}
\noindent in the case $n$ is even.\\\\

\indent\textbf{STEP\,4.}
 We show $\psi_{1}+\psi_{2}+\psi_{3}\simeq*$.\\\\\indent
We choose arbitrary $u\wedge s\wedge a\wedge t\in U\wedge S^{1}\wedge S^{n}\wedge I,\;(S^{1}=[0,1]/\n0,1\nn ),$ for simplicity, we write the symbol  $u\wedge s\wedge a\wedge t$ by $\x s,t\xx$, and write the symbol  $u\wedge t\wedge a\wedge s$ by $\x t, s\xx$.\\\indent
By  formula \ref{p1d}, \ref{p1dp}, \ref{ap3} and \ref{ap3p}, we obtain, in the case $n$ is odd,

$$\psi_{1}\x s,t\xx=\psi_{3}\x t,s\xx
$$
\noindent
if  $0\leq s\leq1/2$, $1/2\leq t\leq1$; in the case $n$ is even,$$\psi_{1}\x s,t\xx=\psi_{3}\x t,s\xx
$$
\noindent
if  $1/2\leq s\leq1$, $1/2\leq t\leq1$.  By formula \ref{gx}, we have
  \[\mathcal{G}'_{1-2t}=\begin{cases}
\mathcal{G}_{1-4t},\;&\text{    $t\in[0,1/4]$}\\
 \con (p_{1}\s^{n+1}f_{4})(-\wedge (1-2t)) ,\;&\text{$t\in[1/4,1/2]$}
\end{cases}\]\[\mathscr{F}'_{2t-1}=\begin{cases}
\mathscr{F}_{4t-3},\;&\text{    $t\in[3/4,1]$}\\
 (\ts^{n}\con(f_{2}\hc(-1)^{n}p_{3}))(-\wedge (2t-1)) ,\;&\text{$t\in[1/2,3/4].$}
\end{cases}\]
\\\noindent

\noindent We notice for any map $f:\mathrm{C}V\rightarrow W$, the map $f^{+\frac{1}{2}}|_{V\wedge [0,\frac{1}{2}]}$ is constant with respect to $t\in[0,\frac{1}{2}]$. Then, we obtain,
\begin{eqnarray}
&
&\notag(\psi_{1}+\psi_{3})+\psi_{2}\\\notag&
\simeq & \fl\n f_{1}\hc\mathscr{F}',\;\mathcal{G}'\hc \mathrm{C}^{\s^{n}(-1)^{n}\s f_{5}}  \nn+\fl\n\;(\mathcal{G}\hc \mathrm{C}^{\s^{n}(-1)^{n}\s f_{5}})^{+\frac{1}{2}},\;(f_{1}\hc\mathscr{F})^{+\frac{1}{2}}\;\nn  \\\notag&
\simeq& *
\end{eqnarray}
we notice $[\s^{n+2}U,X]$ is abelian,
then our proposition is established.

\end{proof}

\end{document}